\documentclass[a4paper,11pt]{amsart}

\usepackage{amsfonts,amsthm,amssymb,amsmath,amscd}
\usepackage{latexsym}
\usepackage{mathrsfs}
\usepackage[all]{xy}

\newcommand{\ZZ}{\mathbb{Z}}
\newcommand{\QQ}{\mathbb{Q}}

\newcommand{\CC}{\mathbb{C}}
\newcommand{\FF}{\mathbb{F}}
\newcommand{\PP}{\mathbb{P}}

\newcommand{\HH}{\mathbb{H}}

\newcommand{\Si}{\mathfrak{S}}

\newcommand{\GL}{\mathrm{GL}}

\newcommand{\G}{\Gamma}
\newcommand{\I}{\mathbf{1}}

\newcommand{\z}{\zeta}
\begin{document}
\title{On Jacobian Kummer surfaces}
\author{Kenji Koike, Yamanashi University} 
\address{Faculty of Education \\
University of Yamanashi \\
Takeda 4-4-37, Kofu \\
Yamanashi 400-8510, Japan\\}
\date{}
\email{kkoike@yamanashi.ac.jp}
\begin{abstract}
We give explicit equations of smooth Jacobian Kummer surfaces of degree $8$ by theta functions. 
As byproducts, we can write down Rosenhain's $80$ hyperpanes and $32$ lines on these Kummer surfaces 
explicitly.
\end{abstract}
\maketitle
\section{Introduction}
A $1$-dimensional complex torus $E_{\tau} = \CC / (\tau \ZZ + \ZZ)$ can be embedded in $\PP^3$ by
fourth order theta functions. Explicitly, the holomorphic map
\[
E_{\tau} \longrightarrow \PP^3, \qquad
z \mapsto [\theta_{00}:\theta_{01}:\theta_{10}:\theta_{11}] (2z, \tau)
= [x_{00}:x_{01}:x_{10}:x_{11}]
\]
is an isomorphism from $E_{\tau}$ to a complete intersection
\begin{align*} 
(\ast 1) \ \begin{cases}
 a_0^2 x_{00}^2 = a_1^2 x_{01}^2 + a_2^2 x_{10}^2 \\
 a_0^2 x_{11}^2 = a_2^2 x_{01}^2 - a_1^2 x_{10}^2 
\end{cases}
\end{align*}
with coefficients 
\[
 a_0^2 = \theta_{00}(0,\tau)^2, \quad a_1^2 = \theta_{01}(0,\tau)^2, \quad a_2^2 = \theta_{10}(0,\tau)^2
\]
(see \cite{Mu}). In the higher dimensional case, defining equations of Abelian varieties are very complicated. 
For example, fourth order theta functions embed principally polarized Abelian surfaces into $\PP^{15}$. 
Flynn showed that the Jacobian variety $\mathrm{Jac}(C)$ of 
a curve $C$ of genus two is defined by $72$ quadrics in $\PP^{15}$, and he gave explicit $72$ equations
in terms of coefficients of the equation of $C$ (\cite{Fl}). On the other hand, 
the Kummer surface $\mathrm{Jac}(C)/ \{ \pm 1\}$ is given as a quartic surface in $\PP^3$, and 
its minimal desingularization $\mathrm{Km}(C)$ is given as a complete intersection of three quadrics in $\PP^5$. 
In \cite{Kl}, Klein gave Kummer surfaces 
\begin{align*}
\begin{cases}
x_1^2 + \cdots + x_6^2 = 0 \\
k_1 x_1^2 + \cdots + k_6 x_6^2 = 0 \\
k_1^2 x_1^2 + \cdots + k_6^2 x_6^2 = 0
\end{cases}
\end{align*}
as singular surfaces of quadratic line complexes. These surfaces are desingularized Kummer surfaces of curves 
\[
 y^2 = (x-k_1) \cdots (x-k_6)
\]
of genus two (see \cite{GH}).
\\ \indent
In this paper, we rewrite Klein's models by Riemann's theta functions. Namely, we study a rational map 
from a principally polarized Abelian surface $X$ to $\PP^5$ given by six odd theta functions of order $4$.
Applying Riemann's theta relations, we obtain defining equations of $\mathrm{Km} (X)$:
\[
 \begin{matrix} 
  \text{(E1)} & A_{10}^2 X_1^2 & & & + A_1^2 X_4^2 - A_2^2 X_5^2 - A_5^2 X_6^2 = 0, \\
  \text{(E2)} & & A_{10}^2 X_2^2 & & + A_3^2 X_4^2 - A_4^2 X_5^2 - A_8^2 X_6^2 = 0, \\
  \text{(E3)} & & & A_{10}^2 X_3^2 & + A_6^2 X_4^2 - A_9^2 X_5^2 - A_7^2 X_6^2 = 0,
 \end{matrix}
\]
where $A_1, \cdots, A_{10}$ are ten even theta constants (Proposition \ref{Prop1}). 
These equations are considered as a two-dimensional
analogue of $(\ast 1)$. Note that coefficients $a_0^2, \ a_1^2, \ a_2^2$ in $(\ast 1)$ are 
modular forms of level $4$. In fact, the graded ring of modular forms of level $4$ is given by
\[
 \bigoplus_{k=0}^{\infty} M_k(\Gamma(4)) = \CC[\theta_{00}^2,\ \theta_{01}^2,\ \theta_{10}^2],
 \qquad \theta_{00}^4 - \theta_{01}^4 - \theta_{10}^4 = 0,
\]
and $(\ast 1)$ gives an elliptic fibration over $\HH / \Gamma(4)$. More precisely, $(\ast 1)$ together with 
\begin{align*} (\ast 2) \begin{cases}
a_1^2 x_{00}^2 + a_2^2 x_{11}^2 = a_0^2 x_{01}^2 \\
a_2^2 x_{00}^2 - a_1^2 x_{11}^2 = a_0^2 x_{10}^2 
\end{cases} \end{align*}
define the elliptic modular surfaces $S(4)$, which is isomorphic to the Fermat quartic surface 
(\cite{BH}, \cite{Mu}, \cite{Sh}). 
In our case, there are $15$ quadratic equations (E1), $\cdots$, (E15) in $X_1, \cdots, X_6$ with 
coefficients $A_1^2, \cdots, A_{10}^2$, defining a fibration of Kummer surfaces over the Siegel modular 
$3$-fold $\mathcal{A}_2(2,4)$ (Theorem \ref{theorem1}). 
These $15$ quadrics are those of rank $4$ in the net spanned by (E1), (E2) and (E3). They are
$4$-terms theta relations, and we can find such relations in classical literature (e.g. \cite{Co} and \cite{Kr}).  
In this paper, we determine singular fibers over 
boundaries of $\mathcal{A}_2(2,4)$ (section 4), and we write down $80$ Rosenhain's hyperplanes that cut out
32 lines on Jacobian Kummer surfaces (Theorem \ref{Th-Rosen}). 
\section{Theta functions and Siegel modular varieties}
\subsection{} 
We denote the Siegel upper half space of degree $g$ by $\Si_g$, the symplectic group $\mathrm{Sp}_{2g}(\ZZ)$ by $\G_g$, 
and Igusa's congruence subgroups by $\G_g(2n,4n)$:
\begin{align*}
 \Si_g &= \{ \Omega \in \mathrm{M}_g(\CC) \ | \ {}^t\Omega = \Omega,\ \mathrm{Im} \medspace \Omega > 0 \}, \\
\G_g &= \{ \gamma \in \GL_{2g}(\ZZ) \ | \ {}^t \gamma J \gamma = J \}, \quad 
J = \begin{bmatrix} 0 & -\I_g \\ \I_g & 0\end{bmatrix}, \\
\G_g(n) &= \{ \gamma \in \G_2 \ | \ \gamma \equiv \I_g \mod n \}, \\
\G_g(2n,4n) &= \{ \begin{bmatrix} A & B \\ C & D \end{bmatrix} \in \G_g(2n) \ | \
\mathrm{diag}(B) \equiv \mathrm{diag}(C) \equiv \I_g \mod 4n \}.
\end{align*}
(Usually the symbol $\G_1(n)$ denotes another group, but here we 
use it for the principal congruence subgroup of level $n$ for $g=1$.)
In the case of $g=1$, we have $\G_1(2,4) / \G_1(4) = \{ \pm 1 \}$. For $g=2$, we have
\[
 \G_2 / \G_2(2) \cong \mathrm{Sp}_{4}(\ZZ / 2\ZZ) \cong \mathrm{S}_6, \quad 
\G_2(2) / \G_2(2,4) \cong (\ZZ / 2\ZZ)^4, \quad \G_2(2,4) / \G_2(4,8) \cong (\ZZ / 2\ZZ)^9.
\]
The group $\G_g$ acts on $\CC^g \times \Si_g$ by
\[
 \begin{bmatrix} A & B \\ C & D \end{bmatrix} \cdot (z, \Omega) = ({}^t (C\Omega + D)^{-1} z, \ (A \Omega + B)(C \Omega + D)^{-1}),
\]
and acts on $(\QQ^g / \ZZ^g)^2$ by
\[
 \begin{bmatrix} A & B \\ C & D \end{bmatrix} \cdot \begin{bmatrix} a \\ b \end{bmatrix} = 
 \begin{bmatrix} D a - C b + \frac{1}{2} \mathrm{diag}(C^tD) \\ -B a + A b + \frac{1}{2} \mathrm{diag}(A^tB) \end{bmatrix},
\quad a, b \in \QQ^g / \ZZ^g.
\]
Theta functions with characteristics $a, b \in \QQ^g$ are defined by
\[
 \vartheta \begin{bmatrix} a \\ b \end{bmatrix}(z, \Omega) 
= \sum _{n \in \ZZ^g} \exp [ \pi i \medspace {}^t(n+a) \Omega (n+a) + 2 \pi i \medspace {}^t(n+a) (z+b)] 
\]
and they satisfy the automorphy property(the theta transformation formula):
\begin{align*}
\vartheta \begin{bmatrix} a^{\sharp} \\ b^{\sharp} \end{bmatrix} (z^{\sharp}, \Omega^{\sharp}) 
&= \kappa(\gamma) \det(C \Omega + D)^{\frac{1}{2}} F(a, b, g, \Omega, z)
\vartheta \begin{bmatrix} a \\ b \end{bmatrix} (z, \Omega)
\end{align*}
where $\gamma = \begin{bmatrix} A & B \\ C & D \end{bmatrix} \in \G_g$, \
$\begin{bmatrix} a^{\sharp} \\ b^{\sharp} \end{bmatrix} = \gamma \cdot \begin{bmatrix} a \\ b \end{bmatrix}$,\ 
$(z^{\sharp}, \Omega^{\sharp}) = \gamma \cdot (z, \Omega)$, $\kappa(\gamma)^8=1$ and
\begin{align*}
 F(a, b, \gamma, \Omega, z) 
= \exp[\pi i \{^t(D a - C b)(-Ba + Ab + (A^tB)_0) - {}^ta b + {}^tz (C \Omega + D)^{-1} C z\}].
\end{align*}
(see \cite{BL} \S 8.5 and \S 8.6).
We can embed Siegel modular 3-folds 
\[
 \mathcal{A}_2(2n,4n) = \mathfrak{S}_2 / \G_2(2n,4n)
\]
for $n=1, 2$ into projective spaces by
theta constants with half integer characteristics. For the simplicity, 
we denote theta constants
\[
 \vartheta \begin{bmatrix} a \\ b \end{bmatrix}(0, \Omega) \quad \text{with} \quad 
a = \frac{1}{2} \begin{bmatrix} i \\ j \end{bmatrix}, \
b = \frac{1}{2} \begin{bmatrix} k \\l \end{bmatrix} \in \frac{1}{2} \ZZ^2
\]
by $\theta [^{ij}_{kl}](\Omega)$. Then we have
\subsection{Proposition (\cite{vGN}, \cite{vGvS}, \cite{Ig2}, \cite{Sa})}
\label{vG}
{\bf (1)} \ The holomorphic map
\[
 \Theta_{2,4} : \Si_2 \longrightarrow \PP^3, \quad \Omega \mapsto
[\theta[_{00}^{00}](2\Omega):\theta[_{00}^{01}](2\Omega):\theta[_{00}^{10}](2\Omega):\theta[_{00}^{11}](2\Omega)]
= [B_0:B_1:B_2:B_3]
\]
gives an isomorphism of the Satake compactification $\overline{\mathcal{A}_2(2,4)}$ and $\PP^3$.
\\
{\bf (2)} \ The holomorphic map
\[
 \Theta_{4,8} : \Si_2 \longrightarrow \PP^9, \quad \Omega \mapsto
[A_1: \cdots :A_{10}]
\]
by $10$ even theta constants 
\begin{align*}
A_1 = \theta[_{00}^{00}](\Omega), \quad A_2 = \theta[^{00}_{01}](\Omega), \quad A_3 = \theta[^{00}_{10}](\Omega), \quad
A_4 = \theta[^{00}_{11}](\Omega), \quad A_5 = \theta[^{01}_{00}](\Omega), \\
A_6 = \theta[^{10}_{00}](\Omega), \quad A_7 = \theta[^{11}_{00}](\Omega), \quad A_8 = \theta[^{01}_{10}](\Omega), \quad
A_9 = \theta[^{10}_{01}](\Omega), \quad A_{10} = \theta[^{11}_{11}](\Omega),
\end{align*}
gives an isomorphism of the Satake compactification $\overline{\mathcal{A}_2(4,8)}$ and 
the closure of the image of the map $\Theta_{4,8}$.
\\
{\bf (3)} \ We have quadratic relations
\begin{align*}
A_1^2 = B_0^2 + B_1^2 + B_2^2 + B_3^2, \quad A_2^2 = B_0^2 - B_1^2 + B_2^2 - B_3^2, \\
A_3^2 = B_0^2 + B_1^2 - B_2^2 - B_3^2, \quad A_4^2 = B_0^2 - B_1^2 - B_2^2 + B_3^2, \\
A_5^2 = 2(B_0 B_1 + B_2 B_3), \quad 
A_6^2 = 2(B_0 B_2 + B_1 B_3), \quad A_7^2 = 2(B_0 B_3 + B_1 B_2), \\
A_8^2 = 2(B_0 B_1 - B_2 B_3), \quad A_9^2 = 2(B_0 B_2 - B_1 B_3), \quad A_{10}^2 = 2(B_0 B_3 - B_1 B_2)
\end{align*}
and the following diagram
\[
\xymatrix{
\overline{\mathcal{A}_2(4,8)} \ar[rr]^{\Theta_{4,8}} \ar[d]_{\pi} & 
 & \PP^9 \ar[d]^{Sq} \\
\overline{\mathcal{A}_2(2,4)} \ar[r]^{\phantom{aa} \Theta_{2,4}} & \PP^3 \ar[r]^{Ver} & \PP^9
}
\]
where $\pi$ is the natural covering map 
with the covering group $\G_2(2,4) / \G_2(4,8) \cong (\ZZ / 2 \ZZ)^9$, the map $Sq$ is the 
squaring of coordinates $[A_1: \cdots : A_{10}] \mapsto [A_1^2 : \cdots : A_{10}^2]$ and 
the map $Ver$ is the Veronese map defined by the above quadratic relations.
\\
{\bf(4)} The $10$ smooth quadrics $Q_i = \{ A_i = 0\}$ of $\PP^3$ correspond to the closure of the locus of 
decomposable principally polarized Abelian surfaces. Therefore $\mathcal{U} = \PP^3 - \cup_{i=1}^{10} Q_i$ 
parametrize Jacobians of curves of genus two. The intersection $Q_i \cap Q_j$ consists of four lines. 
There are $30$ such lines $L_1, \cdots, L_{30}$. They form
the $1$-dimensional boundaries of $\mathcal{A}_2(2,4)$;
\[
 \overline{\mathcal{A}_2(2,4)} - \mathcal{A}_2(2,4) = \bigcup_{i=1}^{30} L_i.
\]
There are $60$ intersecton points of $L_i$'s. These points are $0$-dimensional boundaries.
\section{Quadratic relations of odd theta functions}
\subsection{}
Let $C$ be a curve of genus two, $X = \mathrm{Jac}(C)$ be its Jacobian and $\theta \cong C$ be the theta divisor 
on $X$. Let $\pi : \tilde{X} \rightarrow X$ be the blow up of $2$-torsion points $p_1, \cdots, p_{16}$, and $E_1, \cdots, E_{16}$ 
be the exceptional divisors.
The linear system $|4 \pi^* \theta - \sum E_i|$ gives a morphism of degree $2$ from $\tilde{X}$ to a complete intersection 
$Y$ of three quadrics in $\PP^5$. In fact, the image $Y$ is the desingularized Kummer surface $\tilde{X} / \{ \pm 1\}$ 
(\cite[Chap. 6]{GH}). Let $N \cong \PP^2$ be the net spanned by these quadrics. 
It is known that the discriminant locus $\Delta \subset N$ corresponding to singular quadrics is a union of $6$ lines. 
There are $15$ intersection points of these $6$ lines corresponding to quadrics of rank $4$ (see e.g. \cite[Theorem 3.3]{CR}).
\\ \indent
Now let $X$ be a complex torus $\CC^2 / \Omega  \ZZ^2 + \ZZ^2$, and $\theta$ be given by $\theta [^{00}_{00}](z, \Omega) = 0$. 
Then $\Gamma(X, \mathcal{O}_X(4\theta))$ is identified with the vector space of $4$-th order theta functions, 
and a basis is given by $16$ theta functions $\theta[^{ij}_{kl}](2z, \Omega)$.
The linear system $|4 \theta - \sum p_i|$ corresponds to the subspace spanned by $6$ odd theta functions 
\begin{align*}
 X_1 = \theta [^{01}_{01}](2z, \Omega), \quad  X_2 = \theta [^{01}_{11}](2z, \Omega), \quad X_3 = \theta [^{11}_{01}](2z, \Omega), \\
X_4 = \theta [^{10}_{10}](2z, \Omega), \quad X_5 = \theta [^{10}_{11}](2z, \Omega), \quad X_6 = \theta [^{11}_{10}](2z, \Omega).
\end{align*}
Therefore the defining equations of $Y$ are given by theta relations of $X_1, \cdots, X_6$. 
We can find $4$-terms quadratic relations of theta functions in \cite[\S 30]{Co} and \cite[Chap. VII, \S 14]{Kr}, 
and they will give $15$ quadrics of rank $4$. Here we write down $15$ equations in our terms.
\subsection{}
To get quadratic relations, we apply Riemann's relation \cite[p.214, $\mathrm{(R_{CH})}$]{Mu}
\begin{align*}
\vartheta \begin{bmatrix} \frac{a+b+c+d}{2} \\ \frac{e+f+g+h}{2} \end{bmatrix}(\frac{x+y+u+v}{2}, \Omega)
\cdots
\vartheta \begin{bmatrix} \frac{a-b-c+d}{2} \\ \frac{e-f-g+h}{2} \end{bmatrix}(\frac{x-y-u+v}{2}, \Omega) \\
= \frac{1}{4} \sum_{\alpha, \beta \in \frac{1}{2} \ZZ^2 / \ZZ^2} 
\exp[-2 \pi i \medspace {}^t \beta(a+b+c+d)] 
\vartheta \begin{bmatrix} a + \alpha \\ e + \beta \end{bmatrix}(x, \Omega)
\cdots
\vartheta \begin{bmatrix} d + \alpha \\ h + \beta \end{bmatrix}(v, \Omega).
\end{align*}
Putting $x=y,\ u=v=0,\ a = b + p, \ e = f + q, \ c=d=g=h=0$ with $p, q \in \ZZ^2, \ b, f \in \frac{1}{2}\ZZ^2$, we have
\begin{align*}
&\vartheta \begin{bmatrix} b + \frac{1}{2}p \\ f + \frac{1}{2}q \end{bmatrix}(x, \Omega)^2
\vartheta \begin{bmatrix} \frac{1}{2}p \\ \frac{1}{2}q \end{bmatrix}(0, \Omega)^2  
\\
=&\frac{1}{4} \sum_{\alpha, \beta \in \frac{1}{2} \ZZ^2 / \ZZ^2} 
\exp[-2 \pi i \medspace {}^t \beta (2b + p)] 
\vartheta \begin{bmatrix} b+p + \alpha \\ f+q + \beta \end{bmatrix}(x, \Omega)
\vartheta \begin{bmatrix} b + \alpha \\ f + \beta \end{bmatrix}(x, \Omega)
\vartheta \begin{bmatrix} \alpha \\ \beta \end{bmatrix}(0, \Omega)^2 
\\
=&\frac{1}{4} \sum_{\alpha, \beta \in \frac{1}{2} \ZZ^2 / \ZZ^2} 
(-1)^{2(^t \beta p + ^t\alpha q)} (-1)^{2^tb q}(-1)^{4^t \beta b}
\vartheta \begin{bmatrix} b + \alpha \\ f + \beta \end{bmatrix}(x, \Omega)^2
\vartheta \begin{bmatrix} \alpha \\ \beta \end{bmatrix}(0, \Omega)^2. 
\end{align*}
Now put $b=f = \begin{bmatrix} \frac{1}{2} \\ 0 \end{bmatrix}$ and 
$S = \{\begin{bmatrix} 0 \\ 0 \end{bmatrix}, \ \begin{bmatrix} 0 \\ 1 \end{bmatrix}\}$. We have ${}^t bq = 0$ 
for $q \in S$, and therefore
\begin{align*}
&\sum_{p,q \in S} (-1)^{{}^tpp + {}^tqq }  
\vartheta \begin{bmatrix} b + \frac{1}{2}p \\ f + \frac{1}{2}q \end{bmatrix}(x, \Omega)^2
\vartheta \begin{bmatrix} \frac{1}{2}p \\ \frac{1}{2}q \end{bmatrix}(0, \Omega)^2
\\
&= \sum_{p,q \in S} (-1)^{{}^tpp + {}^tqq } \frac{1}{4} \sum_{\alpha, \beta \in \frac{1}{2} \ZZ^2 / \ZZ^2} 
(-1)^{2(^t \beta p + ^t\alpha q)} (-1)^{4^t \beta b}
\vartheta \begin{bmatrix} b + \alpha \\ f + \beta \end{bmatrix}(x, \Omega)^2
\vartheta \begin{bmatrix} \alpha \\ \beta \end{bmatrix}(0, \Omega)^2
\\
&= \frac{1}{4} \sum_{\alpha, \beta \in \frac{1}{2} \ZZ^2 / \ZZ^2} 
\left( \sum_{p,q \in S} (-1)^{{}^t (2\beta +p)p + {}^t (2\alpha +q)q} \right)
(-1)^{4^t \beta b}
\vartheta \begin{bmatrix} b + \alpha \\ f + \beta \end{bmatrix}(x, \Omega)^2
\vartheta \begin{bmatrix} \alpha \\ \beta \end{bmatrix}(0, \Omega)^2,
\end{align*}
where
\begin{align*}
\sum_{p,q \in S} (-1)^{{}^t (2\beta +p)p + {}^t (2\alpha +q)q} = 
\begin{cases} 4 \quad (\alpha, \beta = (*,\frac{1}{2})) \\
0 \quad (\text{others})\end{cases}.
\end{align*}
Since odd theta constants vanish, the result is 
\begin{align*}
&\theta[^{10}_{10}](x, \Omega)^2 \theta[^{00}_{00}](0, \Omega)^2 - \theta[^{11}_{10}](x, \Omega)^2 \theta[^{01}_{00}](0, \Omega)^2
- \theta[^{10}_{11}](x, \Omega)^2 \theta[^{00}_{01}](0, \Omega)^2 
\\
= &-\theta[^{01}_{01}](x, \Omega)^2 \theta[^{11}_{11}](0, \Omega)^2.
\end{align*}
Putting $x=2z$, we have a quadratic relation 
\[
 \text{(E1)} \quad A_{10}^2 X_1^2 + A_1^2 X_4^2 - A_2^2 X_5^2 - A_5^2 X_6^2 = 0.
\]
Since $\G_2$ acts on the projective coordinates $[A_1: \cdots:A_{10}]$ and $[X_1 : \cdots : X_6]$ via the transformation formula
of theta functions, we can easily derive other $14$ quadrics of rank $4$ from (E1). We summurize the action of $\G_2$ and 
$15$ equations in Appendix. Now we have
\subsection{Proposition}
\label{Prop1}
Let $X = \CC^2 / \Omega \ZZ^2 + \ZZ^2$ be a principally polarized Abelian surface with
$A_i \ne 0$ for $i=1, \cdots, 10$ (hence $X$ is the Jacobian of a curve of genus $2$).  
Then the image $Y$ of a rational map $\Phi : X \longrightarrow \PP^5$
\begin{align*}
 z \mapsto &[ \theta[^{01}_{01}] : \theta [^{01}_{11}] : \theta [^{11}_{01}] :
\theta [^{10}_{10}] : \theta [^{10}_{11}] : \theta [^{11}_{10}]](2z, \Omega) = [X_1:X_2:X_3:X_4:X_5:X_6]
\end{align*}
is a smooth complete intersection given by
\[
 \begin{matrix} 
  \text{(E1)} & A_{10}^2 X_1^2 & & & + A_1^2 X_4^2 - A_2^2 X_5^2 - A_5^2 X_6^2 = 0, \\
  \text{(E2)} & & A_{10}^2 X_2^2 & & + A_3^2 X_4^2 - A_4^2 X_5^2 - A_8^2 X_6^2 = 0, \\
  \text{(E3)} & & & A_{10}^2 X_3^2 & + A_6^2 X_4^2 - A_9^2 X_5^2 - A_7^2 X_6^2 = 0,
 \end{matrix}
\]
and the net spanned by (E1), (E2) and (E3) contains $15$ quadrics of rank $4$ deinfed by (E1), $\cdots$, (E15) 
in Appendix.  
\subsection{Remark}
\label{Rem-sign}
Translations by $2$-torsion points $\frac{1}{2}\Omega \ZZ^2 + \frac{1}{2} \ZZ^2$
acts on the projective coordinate of $\PP^5$ 
as sign changes:
\begin{center} \begin{tabular}{c||c|c|c|c|c|c}
 & $X_1$ & $X_2$ & $X_3$ & $X_4$ & $X_5$ & $X_6$
\\ \hline \hline
$\frac{1}{2} \Omega [^1_0]$ & $+$ & $-$ & $+$ & $-$ & $-$ & $-$  
\\ \hline
$\frac{1}{2} \Omega [^0_1]$ & $-$ & $-$ & $-$ & $+$ & $-$ & $+$
\\ \hline
$\frac{1}{2}[^1_0]$ & $+$ & $+$ & $-$ & $-$ & $-$ & $-$
\\ \hline
$\frac{1}{2}[^0_1]$ & $-$ & $-$ & $-$ & $+$ & $+$ & $-$
\end{tabular}
\end{center}
\subsection{Remark}
\label{Rem-Rosen}
The Resenhain normal form (\cite{Ig1}) of a curve of genus $2$ is given by
\[
 y^2 = x(x-1)(x-\lambda_1)(x-\lambda_2)(x-\lambda_3)
\]
with
\[
 \lambda_1 = \frac{A_7^2 A_6^2}{A_5^2 A_1^2}, \quad \lambda_2 = \frac{A_9^2 A_7^2}{A_2^2 A_5^2}, \quad
\lambda_3 = \frac{A_9^2 A_6^2}{A_2^2 A_1^2}.
\]
\section{Fibration of Kummer surfaces}
\subsection{}
Let us recall that coefficients $A_1^2, \cdots, A_{10}^2$ of $15$ equations (E1), $\cdots$, (E15) are
quadric polynomials in $B_0, B_1, B_2, B_3$ (see Proposition 2.2). Therefore these equations define
a projective variety $\mathcal{X} \subset \PP^3 \times \PP^5$, and we have the projection 
\[
 \pi_{2,4} : \mathcal{X} \longrightarrow \PP^3 \cong \overline{\mathcal{A}_2(2,4)}, \quad
[B_0: \cdots : B_3] \times [X_1: \cdots : X_6] \mapsto [B_0: \cdots : B_3].
\] 
Over $\mathcal{U} = \PP^3 - \cup_{i=1}^{10} Q_i$, this is a K3-fibration and $\pi_{2,4}^{-1}(\mathcal{U})$ is 
smooth. Let us investigate fibers over $\cup_{i=1}^{10}Q_i$. Since $\G_2$ acts on $10$ quadrics 
$\{ Q_1, \cdots, Q_{10} \}$ doubly transitive, we look at fibers only over 
$Q_{10} = \{B_0 B_3 - B_1 B_2 = 0 \}$. We identify $Q_{10}$ with $\PP^1 \times \PP^1$ by the Segre embedding
\[
 Seg : \PP^1 \times \PP^1 \longrightarrow \PP^3, \qquad 
[s_0:s_1] \times [t_0:t_1] \mapsto [s_0 t_0 : s_0 t_1 : s_1 t_0 : s_1 t_1].
\]
Then we have a commutative diagram 
\[
\begin{CD} 
\HH / \G_1(4) \times \HH / \G_1(4) @>\psi>> \PP^1 \times \PP^1 \\
@V\varepsilon VV @VVSegV \\
\Si_2 / \G_2(2,4) @>\Theta_{2,4}>> \PP^3 
\end{CD} 
\]
where $\varepsilon(\tau_1, \tau_2) = \begin{bmatrix} \tau_1 & 0 \\ 0 & \tau_2 \end{bmatrix}$ and $\psi$ is given by
\begin{align*}
(\tau_1, \tau_2) \mapsto [\theta_{00}(2\tau_1) :\theta_{10}(2\tau_1)] \times [\theta_{00}(2\tau_2) :\theta_{10}(2\tau_2)].
\end{align*}
\subsection{Remark}
The duplication formula for $g=1$ are   
\begin{align*}
 \theta_{00}(\tau)^2 = \theta_{00}(2\tau)^2 + \theta_{10}(2\tau)^2, \qquad 
 \theta_{10}(\tau)^2 = 2 \theta_{00}(2\tau) \theta_{10}(2\tau), \\
 \theta_{01}(\tau)^2 = \theta_{00}(2\tau)^2 - \theta_{10}(2\tau)^2.
\end{align*}
\subsection{}
Since we have $\HH / \G_1(4) \cong \PP^1 - \{ \text{6 points}\}$ (\cite{Mu}), the boundary 
$\PP^1 \times \PP^1 - \mathrm{Im} \psi$ is decomposed into $6+6$ lines 
(these $12$ lines are obtained as intersections of $Q_{10}$ and $Q_i$); 
\begin{align*}
[B_0:B_1:B_2:B_3] = \begin{cases} 
[t_0:t_1:0:0] \quad ([s_0:s_1] = [1:0]) \\
[0:0:t_0:t_1] \quad ([s_0:s_1] = [0:1]) \\
[t_0:t_1:\pm t_0:\pm t_1] \quad ([s_0:s_1] = [1:\pm1]) \\
[t_0:t_1:\pm \sqrt{-1} t_0:\pm \sqrt{-1} t_1] \quad ([s_0:s_1] = [1:\pm \sqrt{-1}]) \\
\end{cases} 
\\
[B_0:B_1:B_2:B_3] = \begin{cases} 
[s_0:0:s_1:0] \quad ([t_0:t_1] = [1:0]) \\
[0:s_0:0:s_1] \quad ([t_0:t_1] = [0:1]) \\
[s_0:\pm s_0:s_1:\pm s_1] \quad ([t_0:t_1] = [1:\pm1]) \\
[s_0:\pm \sqrt{-1} s_0:s_1:\pm \sqrt{-1} s_1] \quad ([t_0:t_1] = [1:\pm \sqrt{-1}]) \\
\end{cases}
\end{align*}
and they intersect at $6 \times 6$ points.
\subsection{}
For $\Omega = \varepsilon(\tau_1, \tau_2)$, we have
\[
 X = \CC^2 / (\Omega \ZZ^2 + \ZZ^2) = E_1 \times E_2, \qquad E_i = \CC / (\tau_i \ZZ + \ZZ).
\]
The rational map $\Phi : X \dashrightarrow \PP^5$ is the composition of an embedding 
\begin{align*}
E_1 &\times E_2 \longrightarrow \PP^3 \times \PP^3, 
\\
 &(z_1, z_2) \mapsto [\theta_{00}:\theta_{01}:\theta_{10}:\theta_{11}] (2z_1, \tau_1) \times 
[\theta_{00}:\theta_{01}:\theta_{10}:\theta_{11}] (2z_2, \tau_2)
\end{align*}
and a rational map
\begin{align*}
 \PP^3 &\times \PP^3 \dashrightarrow \PP^5, 
\\
&[x_0:x_1:x_2:x_3] \times [y_0:y_1:y_2:y_3] \mapsto 
[x_0 y_3:x_1 y_3:x_2 y_3:x_3 y_0:x_3 y_1:x_3 y_2].
\end{align*}
The corresponding fiber $\mathcal{X}_{(s,t)}$ over 
$[s_0 t_0 : s_0 t_1 : s_1 t_0 : s_1 t_1] \in Q_{10}-\{12 \ \text{lines}\}$ is defined by
\[
  \begin{cases} 
  T_{00} X_4^2 - T_{01} X_5^2 - T_{10} X_6^2 = 0 \quad (\text{E1}), \\
  S_{00} X_1^2 - S_{01} X_2^2 - S_{10} X_3^2 = 0 \quad (\text{E4}), \\
  T_{10} (S_{01} X_1^2 - S_{00} X_2^2 ) - S_{10} (T_{01} X_4^2 - T_{00} X_5^2) = 0 \quad (\text{E7}) 
 \end{cases}
\]
where 
\begin{align*}
 S_{00} = s_0^2 + s_1^2, \quad S_{10} = 2 s_0 s_1, \quad S_{01} = s_0^2 - s_1^2, \\
T_{00} = t_0^2 + t_1^2, \quad T_{10} = 2 t_0 t_1, \quad T_{01} = t_0^2 - t_1^2
\end{align*}
(compare with Remark 4.2).
It is easily shown that $\mathcal{X}_{(s,t)}$ is birational to the Kummer surface of $E_1 \times E_2$. 
More precisely, $8$ curves
\[
 E_1 \times \{ \text{2-torsions of} \ E_2 \}, \quad 
\{ \text{2-torsions of}\ E_1 \} \times E_2
\]
are corresponding to $8$ nodes
\[
 [X_1: \cdots: X_6] = \begin{cases} [\pm \sqrt{\frac{S_{00}}{S_{10}}}:\pm \sqrt{\frac{S_{01}}{S_{10}}}:1:0:0:0] \\ 
[0:0:0:\pm \sqrt{\frac{T_{00}}{T_{10}}}:\pm \sqrt{\frac{T_{01}}{T_{10}}}:1]\end{cases},
\]
and $4 \times 4$ lines
\[
  \{[\pm z_0 \sqrt{\frac{S_{00}}{S_{10}}}:\pm z_0 \sqrt{\frac{S_{01}}{S_{10}}}:z_0:
\pm z_1 \sqrt{\frac{T_{00}}{T_{10}}}:\pm z_1 \sqrt{\frac{T_{01}}{T_{10}}}:z_1] \ | \ [z_0:z_1] \in \PP^1 \}
\] 
joining $4+4$ nodes are exceptioanl divisors obtained by blowing up at $2$-torsions. 
Resolving $8$ nodes, we obtain the smooth Kummer surface.
\subsection{Remark}
The $8$ curves on $E_1 \times E_2$ are fixed loci by 
\[
 \pm \begin{bmatrix} \alpha & 0 \\ 0 & {}^t\alpha^{-1} \end{bmatrix} \in \G_2(2,4), \qquad
\alpha = \begin{bmatrix} -1 & 0 \\ 0 & 1\end{bmatrix}.
\] 
\subsection{$1$-dimensional boundary}
Next we investugate fibers over a $1$-dimensional boundary
\[
 [B_0:B_1:B_2:B_3] = [t_0:t_1:0:0] \qquad ([s_0:s_1] = [1:0]).
\]
Outside six $0$-dimensional boundaries
\[
 t = t_0 / t_1 = 0, \ \infty, \ \pm1, \ \pm \sqrt{-1}, 
\]
the fiber $\pi_{2,4}^{-1}(t)$ is defined by
\[
\begin{cases} 
 X_1^2 - X_2^2 = 0 \quad (\text{E}4, 5, 6, 7, 8, 9), \\ 
 T_{00} X_4^2 - T_{01} X_5^2 -T_{10} X_6^2 = 0 \quad (\text{E}1, 2), \\ 
 T_{10} X_3^2 - T_{01} X_4^2 + T_{00} X_5^2 = 0 \quad (\text{E}10, 13), \\
 T_{01} X_3^2 + T_{10} X_4^2 - T_{00} X_6^2 = 0 \quad (\text{E}11, 14), \\
 T_{00} X_3^2 + T_{10} X_5^2 - T_{01} X_6^2 = 0 \quad (\text{E}12, 15).
\end{cases}
\] 
((E3) vanish identically). Because there are linear relations,
\[
 (\text{E}10) = -\frac{T_{01}}{T_{00}} (\text{E}1) + \frac{T_{10}}{T_{00}} (\text{E}12), \qquad 
(\text{E}11) = \frac{T_{10}}{T_{00}} (\text{E}1) + \frac{T_{01}}{T_{00}} (\text{E}12),
\]
we have $\pi_{2,4}^{-1}(t) = Y_+ \cup Y_-$ with 
\[
 Y_+ = \begin{cases} 
X_1 = X_2 \\ T_{00} X_4^2 - T_{01} X_5^2 -T_{10} X_6^2 = 0 \\ T_{00} X_3^2 + T_{10} X_5^2 - T_{01} X_6^2 = 0
\end{cases}, \qquad
 Y_- = \begin{cases} 
X_1 = -X_2 \\ T_{00} X_4^2 - T_{01} X_5^2 -T_{10} X_6^2 = 0 \\ T_{00} X_3^2 + T_{10} X_5^2 - T_{01} X_6^2 = 0
\end{cases}.
\]
These surfaces are cones over an elliptic curve $E = Y_+ \cap Y_-$, and singular at $[1:\pm1:0:0:0:0]$.
\subsection{$0$-dimensional boundary}
Let us investigate the fiber at a $0$-dimensional boundary
\[
P_0 =  [B_0:B_1:B_2:B_3] = [1:0:0:0]
\]
corresponding to $[s_0:s_1] \times [t_0:t_1] = [1:0] \times [1:0]$. Then $15$ equations are 
\[
 \begin{cases} 
   X_1^2 - X_2^2 = 0 \quad (\text{E4, 5, 8, 9}), \\ 
   X_3^2 - X_6^2 = 0 \quad (\text{E11, 12, 14 , 15}), \\ 
   X_4^2 - X_5^2 = 0 \quad (\text{E1, 2, 10, 13})
 \end{cases}
\] 
((E3), (E6) and (E7) vanish identically). The fiber $\pi_{2,4}^{-1}(P_0)$ is a union of $8$ projective plane:
\[
 \PP_{+++} = \begin{cases} X_1 = X_2 \\ X_3 = X_6 \\ X_4 = X_5 \end{cases}, \quad 
 \PP_{++-} = \begin{cases} X_1 = X_2 \\ X_3 = X_6 \\ X_4 = -X_5 \end{cases}, \cdots, \quad
 \PP_{---} = \begin{cases} X_1 = -X_2 \\ X_3 = -X_6 \\ X_4 = -X_5 \end{cases},
\]
and the dual graph (see \cite{Pe} ) is a cube (vertices, edges and faces represent irreducible components, double lines, 
poins of order $4$ respectively).
\begin{figure}[htbp] \begin{center}
\setlength{\unitlength}{1mm}
\begin{picture}(50,45)
\thicklines
\put(10,0){\framebox(35,30){}}
\put(25,15){\framebox(35,30){}}
\put(10,0){\line(1,1){15}}
\put(45,0){\line(1,1){15}}
\put(10,30){\line(1,1){15}}
\put(45,30){\line(1,1){15}}
\put(10,30){\circle*{1.5}} \put(0,30){$\PP_{+-+}$} 
\put(45,30){\circle*{1.5}} \put(47,28){$\PP_{+--}$}
\put(10,0){\circle*{1.5}} \put(0,0){$\PP_{--+}$} 
\put(45,0){\circle*{1.5}} \put(48,0){$\PP_{---}$} 
\put(25,45){\circle*{1.5}} \put(15,45){$\PP_{+++}$} 
\put(60,45){\circle*{1.5}} \put(61,45){$\PP_{++-}$}
\put(60,15){\circle*{1.5}} \put(61,15){$\PP_{-+-}$}
\put(25,15){\circle*{1.5}} \put(26,16){$\PP_{-++}$}
\end{picture} \end{center} \caption{} \label{Fgcycle}
\end{figure}
\subsection{Theorem}
\label{theorem1}
{\rm (1)} \ The $5$-dimensional projective variety $\mathcal{X} \subset \PP^3 \times \PP^5$ is smooth and simply connected. 
\\
{\rm (2)} \ The projection $\pi_{2,4} : \mathcal{X} \rightarrow \PP^3$ is the half-anticanonical map, and hence
$\mathcal{X}$ has the Kodaira dimension $\kappa(\mathcal{X}) = - \infty$. 
\\ \\ \noindent
{\bf Proof.}
{\rm (1)} \ The smoothness is shown by the Jacobian criterion. Namely, we show that not all of $3$-minors of 
the projective Jacobian matrix 
\[
 \mathcal{J} = \frac{\partial (\text{E1}, \cdots, \text{E15}) }
{\partial(B_0, \cdots, B_3, X_1, \cdots, X_6)}
\]
of the equations (E1), $\cdots$, (E15) vanish.
\\ \indent
We first consider fibers over the $1$-dimensional boundary
\[
 L = \{[t_0:t_1:0:0] \ | \ [t_0:t_1] \in \PP^1,\ T_{00} T_{10} T_{01} \ne 0 \}
\]
Since we have
\[
\left| \frac{\partial (\text{E1, E4, E11}) }{\partial(B_2, B_3,X_1)} \right| = 8 T_{00} T_{01} X_1^5 
\]
at points of $\pi_{2,4}^{-1}(L)$, the variety $\mathcal{X}$ is smooth at points of $\{X_1 \ne 0 \} \cap \pi_{2,4}^{-1}(L)$.
On the other hand, we have
\begin{align*}
 \left| \frac{\partial (\text{E1, E4, E11}) }{\partial(B_0, B_2, X_4)} \right| &=
-8 X_4 (t_0 X_3^2 + t_1 X_4^2)\{-t_1 T_{01} X_4^2 + t_0(T_{00} X_3^2 + T_{10} X_5^2 - T_{01} X_6^2)\} \\
&= 8 t_1 T_{01} X_4^3 (t_0 X_3^2 + t_1 X_4^2), \\
\left| \frac{\partial (\text{E1, E4, E11}) }{\partial(B_1, B_3, X_6)} \right| &=
8 X_6 (t_0 X_4^2 - t_1 X_3^2)\{-t_0 T_{01}X_6^2 + t_1(T_{10} X_3^2 - T_{01} X_4^2 + T_{00} X_5^2)\} \\
&= -8 t_0 T_{01} X_6^3 (t_0 X_4^2 - t_1 X_3^2).
\end{align*} 
at points of $\{X_1 = 0 \} \cap \pi_{2,4}^{-1}(L)$. We see that these values do not vanish at the same time on 
$\{X_1 = 0 \} \cap \pi_{2,4}^{-1}(L)$ (see {\bf 4.6}), and $\mathcal{X}$ is smooth along $\pi_{2,4}^{-1}(L)$.
\\ \indent
Similarly, we have the following $3$-minors
\begin{align*}
 \left| \frac{\partial (\text{E1, E4, E11}) }{\partial(B_2, B_3,X_1)} \right| = 8X_1^5, \quad
 \left| \frac{\partial (\text{E1, E4, E11}) }{\partial(B_1, B_2, X_3)} \right| = -8X_3^5, \quad
 \left| \frac{\partial (\text{E1, E4, E11}) }{\partial(B_1, B_3, X_4)} \right| = 8X_4^5
\end{align*}
at points of $\pi_{2,4}^{-1}(P_0)$, that is, at points
\[
 [X_1:\pm X_1:X_3:X_4:\pm X_4:\pm X_3] \times [1:0:0:0] \in \mathcal{X} \subset \PP^5 \times \PP^3.
\] 
Therefore the Jacobian matrix $\mathcal{J}$ has rank $3$ at points of fibers over 
$0$-dimensional and $1$-dimensional boundaries, and $\mathcal{X}$ is smooth there. 
\\ \indent
Now let $Sing \mathcal{X}$ be the set of singular points of $\mathcal{X}$. Then $\pi_{2,4}(Sing \mathcal{X}) \subset \PP^3$ is
closed subvariety since $\pi_{2,4}$ is proper. As we have seen, it holds 
\[ 
 \pi_{2,4}(Sing \mathcal{X}) \cap Q_{10} \ \subset \ Q_{10} - \{\text{boundaries}\} \cong 
(\HH / \Gamma_1(4)) \times (\HH / \Gamma_1(4))
\]
and $\pi_{2,4}(Sing \mathcal{X}) \cap Q_{10}$ must be isolated points. However, the fibration 
over $(\HH / \Gamma_1(4)) \times (\HH / \Gamma_1(4))$ is locally trivial as a topological space, and
we see that $\pi_{2,4}(Sing \mathcal{X}) \cap Q_{10} = \phi$. By the $\G_2$-symmetry, we have $\pi_{2,4}(Sing \mathcal{X}) \cap Q_{10} = \phi$
for $i =1, \cdots, 10$, and $\mathcal{X}$ is smooth.
\\ \indent
Since we have an exact sequence 
\[
 \pi_1(\text{a general fiber}) \longrightarrow \pi_1(\mathcal{X}) 
\longrightarrow \pi_1(\PP^3) \longrightarrow 1
\]
of fundamental groups (\cite{No}, Lemma 1.5) , $\mathcal{X}$ is simply connected.
\\ \indent
{\rm (2)} \ Let $\mathcal{I}$ be the ideal sheaf of $\mathcal{X} \subset \PP^3 \times \PP^5$. By the adjunction formula, we have
\[
 \omega_{\mathcal{X}} \cong \omega_{\PP^3 \times \PP^5} \otimes \wedge^3 (\mathcal{I} / \mathcal{I}^2)^* 
\] 
where $p_1$ (resp. $p_2$) is the projection to $\PP^3$ (resp. $\PP^5$).
Hence it is enough to show that 
\[
 \wedge^3 (\mathcal{I} / \mathcal{I}^2) \cong p_1^* \mathcal{O}_{\PP^3}(-2) \otimes p_2^* \mathcal{O}_{\PP^5}(-6).
\]
Now the problem is to show equivalence of two divisors on a smooth variety, we may ignore subvarieties of codimension $2$.  
So we restrict ourself to an open set $\pi_{2,4}^{-1}(U_9 \cup U_{10})$ where $U_i = \PP^3 - Q_i$. 
Note that (E1), (E2) and (E3) generate $\mathcal{I}$ over $\pi_{2,4}^{-1}(U_{10})$ since they are given by
\[
 A_{10}^2 \begin{bmatrix} X_1^2 \\ X_2^2 \\ X_3^2 \end{bmatrix} + M \begin{bmatrix} X_4^2 \\ X_5^2 \\ X_6^2 \end{bmatrix} 
= \begin{bmatrix} 0 \\ 0 \\ 0 \end{bmatrix}, \qquad 
M =  \begin{bmatrix} A_1^2 & -A_2^2 & -A_5^2 \\ A_3^2 & -A_4^2 & -A_8^2 \\ A_6^2 & -A_9^2 & -A_7^2 \end{bmatrix}
\] 
and $\det M = A_{10}^6$ as polynomials of $B_0, \cdots, B_3$. 
Similarly, (E5), (E7) and (E8) generate $\mathcal{I}$ over $\pi_{2,4}^{-1}(U_9)$ (see the action of $g_4$ in Appendix). 
These basis are connected on $\pi_{2,4}^{-1}(U_9 \cap U_{10})$ by
\[
 \begin{bmatrix} \text{(E5)} \\ \text{(E7)} \\ \text{(E8)} \end{bmatrix} = \frac{1}{A_{10}^2}
\begin{bmatrix} A_2^2 & -A_4^2 & -A_9^2 \\ A_8^2 & -A_5^2 & 0 \\ A_3^2 & -A_1^2 & 0 \end{bmatrix}
\begin{bmatrix} \text{(E1)} \\ \text{(E2)} \\ \text{(E3)} \end{bmatrix}, \qquad
\det \begin{bmatrix} A_2^2 & -A_4^2 & -A_9^2 \\ A_8^2 & -A_5^2 & 0 \\ A_3^2 & -A_1^2 & 0\end{bmatrix}
= - 8 A_{10}^2 A_9^4, 
\]
namely,
\[
 \text{(E5)} \wedge \text{(E7)} \wedge \text{(E8)} = -8 \frac{A_9^4}{A_{10}^4} 
\text{(E1)} \wedge \text{(E2)} \wedge \text{(E3)}.
\]
Taking a standard affine open cover $V_{i,j} = \{ B_i \ne 0, \ X_j \ne 0 \}$ of $\PP^3 \times \PP^5$,  
and considering coordinates changes for open sets $V_{i,j} \cap \pi_{2,4}^{-1}(U_9)$ and $V_{k,l} \cap \pi_{2,4}^{-1}(U_{10})$, 
we see that 
\[
 \wedge^3 (\mathcal{I} / \mathcal{I}^2) \cong p_1^* \mathcal{O}_{\PP^3}(-2) \otimes p_2^* \mathcal{O}_{\PP^5}(-6).
\]
\hfill $\Box$
\section{The $80$ Rosenhain hyperplanes}
\subsection{}
Let $X$ and $Y$ be as in Proposition \ref{Prop1}. As classically known,  Kummer surface $Y$ contains two families 
of disjoint $16$ lines in $\PP^5$. One is $16$ exceptional curves obtained by blowing up $16$ nodes, 
and another is $16$ tropes, that are the images of translated theta divisors by $2$-torsions
(that is, tropes are given by $\vartheta[^a_b](z,\Omega) =0$ with $a, b \in \frac{1}{2}\ZZ^2$).
Each line intersects with $6$ lines in the opposite family. 
These lines are cut out by special hyperplanes called Rosenhain's hyperplanes. They cut out $8$ lines 
consisting of $4$ exceptional curves and $4$ tropes, and there are $80$ such hyperplanes (see \cite{GH} Chap.6, \cite{Hu} Chap. VII). 
We can write down them by Riemann's theta relations. 
Putting $x=2z, \ y=u=v=0$ in \cite[p.214, $\mathrm{(R_{CH})}$]{Mu}, we have
\begin{align*}
&\vartheta \begin{bmatrix} \frac{a+b+c+d}{2} \\ \frac{e+f+g+h}{2} \end{bmatrix}(z)
\vartheta \begin{bmatrix} \frac{a+b-c-d}{2} \\ \frac{e+f-g-h}{2} \end{bmatrix}(z)
\vartheta \begin{bmatrix} \frac{a-b+c-d}{2} \\ \frac{e-f+g-h}{2} \end{bmatrix}(z)
\vartheta \begin{bmatrix} \frac{a-b-c+d}{2} \\ \frac{e-f-g+h}{2} \end{bmatrix}(z) \\
= &\frac{1}{4} \sum_{\alpha, \beta \in \frac{1}{2} \ZZ^2 / \ZZ^2} 
\exp[-2 \pi i \medspace {}^t \beta(a+b+c+d)] 
\vartheta [^{a + \alpha}_{e + \beta}](2z)
\vartheta [^{b + \alpha}_{f + \beta}](0)
\vartheta [^{c + \alpha}_{g + \beta}](0)
\vartheta [^{d + \alpha}_{h + \beta}](0).
\end{align*}
If the left hand side is a $4$-th order odd theta function, then the right hand side must be a linear combination of 
$X_1, \cdots, X_6$. If this is the case, the above equation represents a hyperplane in $\PP^5$ cutting $4$ tropes.
\subsection{}
For example, let us consider four functions
\[
 f_1(z) = \theta[^{00}_{00}](z,\Omega), \quad f_2(z) = \theta[^{00}_{01}](z,\Omega), \quad
 f_3(z) = \theta[^{01}_{00}](z,\Omega), \quad f_4(z) = \theta[^{01}_{01}](z,\Omega).
\]
The product $F(z) = f_1(z) \cdots f_4(z)$ has the same periodicity with $\theta[^{00}_{00}](z,\Omega)^4$,
and it satisfies $F(-z) = -F(z)$. Namely, $F(z)$ is a $4$-th order odd theta function. In fact, we have
\[
 F(z) = -\frac{1}{4}(A_1 A_2 A_5 X_1 + A_3 A_4 A_8 X_2 + A_6 A_7 A_9 X_3)
\]
by putting 
\[
 a = b = e = g = \frac{1}{2} \begin{bmatrix} 0 \\ 0 \end{bmatrix}, \quad 
c = d = f = h = \frac{1}{2} \begin{bmatrix} 0 \\ 1 \end{bmatrix}
\]
in the above theta relation. Denoting $2$-torsion points $\frac{1}{2} \Omega [^i_j] + \frac{1}{2}[^k_l] \in X$ by $^{ij}_{kl}$, 
we have the following table for the vanishing property of $F$.
\\
\begin{center}
\begin{tabular}{c||c|c|c|c|c|c|c|c|}
 & $^{00}_{00}$ & $^{00}_{01}$ & $^{00}_{10}$ & $^{00}_{11}$ & $^{01}_{00}$ & $^{10}_{00}$ & $^{11}_{00}$ & $^{01}_{10}$  
\\ \hline 
$\theta[^{00}_{00}](z)$ & & & & & & & &
\\ \hline
$\theta[^{00}_{01}](z)$ & & & & & $\bullet$ & & $\bullet$ & $\bullet$
\\ \hline
$\theta[^{01}_{00}](z)$ & & $\bullet$ & & $\bullet$ & & & &
\\ \hline
$\theta[^{01}_{01}](z)$ & $\bullet$ & & $\bullet$ & & & $\bullet$ & &
\\ \hline
\end{tabular}
\vskip0.5cm
\begin{tabular}{c||c|c|c|c|c|c|c|c|}
 & $^{10}_{01}$ & $^{11}_{11}$ & $^{01}_{01}$ & $^{01}_{11}$ & $^{11}_{01}$ & $^{10}_{10}$ & $^{10}_{11}$ & $^{11}_{10}$  
\\ \hline 
$\theta[^{00}_{00}](z)$ & & & $\bullet$ & $\bullet$ & $\bullet$ & $\bullet$ & $\bullet$ & $\bullet$
\\ \hline
$\theta[^{00}_{01}](z)$ & & $\bullet$ & & & & $\bullet$ & $\bullet$ & 
\\ \hline
$\theta[^{01}_{00}](z)$ & $\bullet$ & $\bullet$ & & & & $\bullet$ & & $\bullet$
\\ \hline
$\theta[^{01}_{01}](z)$ & & $\bullet$ & & & & & $\bullet$ & $\bullet$
\\ \hline
\end{tabular}
\end{center}
Namely, $F$ vanishes to order $3$ at $^{11}_{11}$, $^{10}_{10}$, $^{10}_{11}$ and $^{11}_{10}$, 
and cuts out $4$ exceptional curves corresponding to these points.
\subsection{}
In general, a product of four theta function with characteristics in $\frac{1}{2}\ZZ^2 / \ZZ^2$
\[
 \vartheta[^{a'}_{a''}](z,\Omega) \vartheta[^{b'}_{b''}](z,\Omega) \vartheta[^{c'}_{c''}](z,\Omega) \vartheta[^{d'}_{d''}](z,\Omega)
\]
has the same periodicity with $\theta[^{00}_{00}](z)^4$ iff 
\[
 a' + b' + c' + d', \ a'' + b'' + c'' + d'' \in \ZZ,  
\]
and it is an odd function iff
\[
 2(a' \cdot a'' + b' \cdot b'' + c' \cdot c'' + d' \cdot d'') \notin \ZZ.  
\]
There are only $80$ such combinations, and we can find them immediately by computer. 
To state the result, we introduce a few notations. We number characteristics from $1$ to $16$:
\begin{center}
\begin{tabular}{|c|c|c|c|c|c|c|c|c|c|c|c|c|c|c|c|}
\hline
1 & 2 & 3 & 4 & 5 & 6 & 7 & 8 & 9 & 10 & 11 & 12 & 13 & 14 & 15 & 16
\\ \hline 
$^{00}_{00}$ & $^{00}_{01}$ & $^{00}_{10}$ & $^{00}_{11}$ & $^{01}_{00}$ & $^{10}_{00}$ & $^{11}_{00}$ & $^{01}_{10}$ 
& $^{10}_{01}$ & $^{11}_{11}$ & $^{01}_{01}$ & $^{01}_{11}$ & $^{11}_{01}$ & $^{10}_{10}$ & $^{10}_{11}$ & $^{11}_{10}$
\\ \hline
\end{tabular}
\end{center}
We denote exceptional curves by $E_1, \cdots, E_{16}$, and tropes by $D_1, \cdots, D_{16}$ according to 
this numbering. Finally, we write $A_{i,j,k}$ instead of $A_i A_j A_k$ and we denote divisors 
$D_i + \cdots + D_j + E_k + \cdots + E_l$ by $D_{i, \cdots, j} + E_{k, \cdots,l}$. 
\subsection{Theorem}
\label{Th-Rosen}
Let $X$ and $Y$ be as in Proposition \ref{Prop1}, and $X[2]$ be the $2$-torsion subgroup of $X$. \\
(1) Rosenhain's hyperplanes for $Y$ are given by $3$-terms linear equation in $X_1, \cdots, X_6$ with coefficients $A_{i,j,k}$. 
Therefore each hyperplane is invariant under the action of a subgroup of $X[2]$ of order $4$ 
(see Remark \ref{Rem-sign}), and each $X[2]$-orbit contains $4$ hyperplanes. Representatives from twenty orbits are given 
explicitly in Appendix.
\\
(2) The trope $D_1$ is the intersection of four hyperplanes
\begin{align*}
H1 :& \quad A_{1,3,6} X_4 + A_{2,4,9} X_5 + A_{5,7,8} X_6 =0, \\
H2 :& \quad A_{1,10,3} X_3 + A_{5,8,9} X_5 - A_{2,4,7} X_6 =0, \\
H5 :& \quad A_{1,6,10} X_2 + A_{4,5,7} X_5 - A_{2,8,9} X_6 =0, \\
H14 :& \quad A_{5,6,9} X_1 + A_{1,2,7} X_3 - A_{3,4,10} X_6 =0,
\end{align*}
and the exceptional curve $E_{11}$ is the intersection of four hyperplanes
\begin{align*}
H1 :& \quad A_{1,3,6} X_4 + A_{2,4,9} X_5 + A_{5,7,8} X_6 =0, \\
H2 :& \quad A_{1,10,3} X_3 + A_{5,8,9} X_5 - A_{2,4,7} X_6 =0, \\
H5 :& \quad A_{1,6,10} X_2 + A_{4,5,7} X_5 - A_{2,8,9} X_6 =0, \\
H11 :& \quad A_{3,6,10} X_1 + A_{2,7,8} X_5 - A_{4,5,9} X_6 =0.
\end{align*}
Note that other topes and exceptional curves are given by $16$ translations.
\subsection{}
We have studied Kummer surfaces over complex numbers until now. However, our family is defined over $\ZZ$, 
and the result is applied for a field $k$ of $\mathrm{char}(k) \ne 2$. 
In fact, we can easily show that (E1), (E2), (E3) define a smooth complete intersection if
$A_i \ne 0$ for $i=1, \cdots, 10$. It is a interesting problem to construct Kummer surfaces with $32$ lines over small 
fields. For example, Kuwata and Shioda asked the probelm finding all elliptic fibrations on given K3 surfaces in \cite{KS},
and they proposed to find $32$ lines on Kummer surfaces to attack this problem in the case of Kummer surfaces. (The problem for 
Jacobian Kummer surfaces were solved by Kumar in \cite{Ku} over algebraically closed fields of characteristic $0$. )
\\ \indent
Now we can construct a Kummer surface with $32$ lines defined over $k$ if we have 
\[
 [B_0:\cdots:B_3] \times [A_1:\cdots:A_{10}] \in \PP^3(k) \times \PP^9(k)
\]
satisfying $A_i \ne 0 \ (i=1, \cdots, 10)$ and the quadric relations in Proposition \ref{vG}. The author does not know 
whether such a point exists for $k = \QQ$. For finite prime fields $\FF_p$, we do not have such a point 
if  $p = 3,5,7,11,13,17$. 
\subsection{Example}
Let us consider a finite prime field $\FF_{19}$, and $b=[1:3:3:3] \in \PP^3(\FF_{19})$. 
The image of $b$ by the Veronese map $Ver$ in Proposition \ref{vG} is
\begin{align*}
 &[9:11: 11: 11: 5: 5: 5: 7: 7: 7] \\
= &[4:7: 7: 7: -2: -2: -2: 1: 1: 1] \\
= &[2^2:8^2: 8^2: 8^2: 6^2: 6^2: 6^2: 1: 1: 1] \in \PP^9(\FF_{19}),
\end{align*}
and we have a smooth Kummer surafce $Y_b$ over $\FF_{19}$ defined by
\[
 \begin{matrix} 
   X_1^2 & & & + 4 X_4^2 - 7 X_5^2 + 2 X_6^2 = 0, \\
   & X_2^2 & & + 7 X_4^2 - 7 X_5^2 - \phantom{a} X_6^2 = 0, \\
   & & X_3^2 & - 2 X_4^2 - \phantom{a} X_5^2 + 2 X_6^2 = 0
 \end{matrix}
\]
with 
\begin{align*}
D_1 : \begin{cases}
X_4 + 7 X_5 -2 X_6 =0, \\
3 X_3 -6 X_5 + 4 X_6 =0, \\
7 X_2 -3 X_5 + 8 X_6 =0, \\
2 X_1 - X_3 + 7 X_6 =0, 
\end{cases}
\quad
E_{11} : \begin{cases}
X_4 + 7 X_5 -2 X_6 =0, \\
3 X_3 -6 X_5 + 4 X_6 =0, \\
7 X_2 -3 X_5 + 8 X_6 =0, \\
X_1 + X_3 - X_6 = 0.
\end{cases}
\end{align*}
The Rosenhain normal form (Remark \ref{Rem-Rosen}) of the corresponding curve of genus two is 
\[
 y^2 = x(x-1)(x-4)(x-9)(x-11).
\]
\appendix
\section{}
\subsection{}
We look at the action of $\G_2$ on the projective coordinates of $\mathcal{A}_2(4,8)$ and $Y = \mathrm{Km}(X)$. 
For unimodular transformations
\[
 g_i = \begin{bmatrix} \alpha_i & 0 \\ 0 & {}^t\alpha_i^{-1}\end{bmatrix}, \qquad
\alpha_1 = [^{01}_{10}], \ \alpha_2 = [^{10}_{11}], \ \alpha_3 = [^{11}_{01}], \ \alpha_4 = [^{11}_{10}], \
\alpha_5 = [^{01}_{11}], 
\]
translations 
\[
 h_i = \begin{bmatrix} \I & \beta_i \\ 0 & \I \end{bmatrix}, \qquad 
\beta_1 = [^{10}_{00}], \ \beta_2 = [^{01}_{10}], \ \beta_3 = [^{00}_{01}] 
\]
and $J = \begin{bmatrix} 0 & -\I \\ \I & 0\end{bmatrix}$, we have the following table, where $\zeta = \exp(2 \pi i /8)$.
\begin{center}
\begin{tabular}{c||c|c|c|c|c|c|c|c|c|c}
$1$ & $A_1$ & $A_2$ & $A_3$ & $A_4$ & $A_5$ & $A_6$ & $A_7$ & $A_8$ & $A_9$ & $A_{10}$  
\\ \hline \hline
$g_1$ & $A_1$ & $A_3$ & $A_2$ & $A_4$ & $A_6$ & $A_5$ & $A_7$ & $A_9$ & $A_8$ & $A_{10}$
\\ \hline
$g_2$ & $A_1$ & $A_2$ & $A_4$ & $A_3$ & $A_7$ & $A_6$ & $A_5$ & $A_{10}$ & $A_9$ & $-A_{8}$
\\ \hline
$g_3$ & $A_1$ & $A_4$ & $A_3$ & $A_2$ & $A_5$ & $A_7$ & $A_6$ & $A_{8}$ & $A_{10}$ & $-A_{9}$
\\ \hline
$g_4$ & $A_1$ & $A_3$ & $A_4$ & $A_2$ & $A_7$ & $A_5$ & $A_6$ & $A_{10}$ & $A_{8}$ & $-A_{9}$
\\ \hline
$g_5$ & $A_1$ & $A_4$ & $A_2$ & $A_3$ & $A_6$ & $A_7$ & $A_5$ & $A_{9}$ & $A_{10}$ & $-A_{8}$
\\ \hline \hline
$h_1$ & $A_3$ & $A_4$ & $A_1$ & $A_2$ & $A_8$ & $\z^7 A_6$ & $\z^7 A_7$ & $A_5$ & $\z^7 A_9$ & $\z^7 A_{10}$
\\ \hline
$h_2$ & $A_1$ & $A_2$ & $A_3$ & $A_4$ & $A_8$ & $A_9$ & $iA_{10}$ & $A_5$ & $A_6$ & $iA_{7}$
\\ \hline
$h_3$ & $A_2$ & $A_1$ & $A_4$ & $A_3$ & $\z^7 A_5$ & $A_9$ & $\z^7 A_7$ & $\z^7 A_8$ & $A_6$ & $\z^7 A_{10}$
\\ \hline \hline
$J$ & $A_1$ & $A_5$ & $A_6$ & $A_7$ & $A_2$ & $A_3$ & $A_4$ & $A_9$ & $A_8$ & $-A_{10}$  
\end{tabular}
\vskip0.5cm
\begin{tabular}{c||c|c|c|c|c|c}
$1$ & $X_1$ & $X_2$ & $X_3$ & $X_4$ & $X_5$ & $X_6$ 
\\ \hline \hline
$g_1$ & $X_4$ & $X_5$ & $X_6$ & $X_1$ & $X_2$ & $X_3$  
\\ \hline
$g_2$ & $X_3$ & $-X_6$ & $X_1$ & $X_5$ & $X_4$ & $X_2$
\\ \hline
$g_3$ & $X_2$ & $X_1$ & $X_5$ & $X_6$ & $-X_3$ & $X_4$ 
\\ \hline
$g_4$ & $X_6$ & $-X_3$ & $X_4$ & $X_2$ & $X_1$ & $X_5$ 
\\ \hline
$g_5$ & $X_5$ & $X_4$ & $X_2$ & $X_3$ & $-X_6$ & $X_1$  
\\ \hline \hline
$h_1$ & $X_2$ & $X_1$ & $\z^7 X_3$ & $\z^7 X_4$ & $\z^7 X_5$ & $\z^7 X_6$ 
\\ \hline
$h_2$ & $X_2$ & $X_1$ & $-iX_6$ & $X_5$ & $X_4$ & $-iX_3$  
\\ \hline 
$h_3$ & $\z^7 X_1$ & $\z^7 X_2$ & $\z^7 X_3$ & $X_5$ & $X_4$ & $\z^7 X_6$ 
\\ \hline \hline
$J$ & $iX_1$ & $iX_3$ & $iX_2$ & $iX_4$ & $iX_6$ & $iX_5$
\end{tabular}
\end{center}
\subsection{}
From the above table, we see that the equation (E1) is transformed in the following equations: 
\begin{align*}
\text{(E1)} \quad A_{10}^2 X_1^2 + A_1^2 X_4^2 - A_2^2 X_5^2 - A_5^2 X_6^2 = 0 \\
\text{(E2)} =h_1 \text{(E1)} \quad A_{10}^2 X_2^2 + A_3^2 X_4^2 - A_4^2 X_5^2 - A_8^2 X_6^2 = 0 \\
\text{(E3)} = J \text{(E2)} \quad A_{10}^2 X_3^2 + A_6^2 X_4^2 - A_9^2 X_5^2 - A_7^2 X_6^2 = 0 \\
\text{(E4)} = g_1 \text{(E1)} \quad A_1^2 X_1^2 - A_3^2 X_2^2 - A_6^2 X_3^2  + A_{10}^2 X_4^2 = 0 \\
\text{(E5)} = g_1 \text{(E2)} \quad A_2^2 X_1^2 - A_4^2 X_2^2 - A_9^2 X_3^2  + A_{10}^2 X_5^2 = 0 \\
\text{(E6)} = g_1 \text{(E3)} \quad A_5^2 X_1^2 - A_8^2 X_2^2 - A_7^2 X_3^2  + A_{10}^2 X_6^2 = 0 \\
\text{(E7)} = g_4 \text{(E3)} \quad A_8^2 X_1^2 - A_5^2 X_2^2 - A_9^2 X_4^2 + A_{6}^2 X_5^2 = 0 \\
\text{(E8)} = g_4 \text{(E1)} \quad A_3^2 X_1^2 - A_1^2 X_2^2 + A_7^2 X_5^2 - A_{9}^2 X_6^2 = 0 \\
\text{(E9)} = h_3 \text{(E8)} \quad A_4^2 X_1^2 - A_2^2 X_2^2 + A_7^2 X_4^2 - A_6^2 X_6^2 = 0\\
\text{(E10)} = J \text{(E9)} \quad A_7^2 X_1^2 - A_5^2 X_3^2 + A_4^2 X_4^2 - A_3^2 X_5^2 = 0\\
\text{(E11)} = g_5 \text{(E2)} \quad A_9^2 X_1^2 - A_2^2 X_3^2 - A_{8}^2 X_4^2 + A_3^2 X_6^2  = 0 \\
\text{(E12)} = g_5 \text{(E1)} \quad A_6^2 X_1^2 - A_1^2 X_3^2 - A_{8}^2 X_5^2 + A_4^2 X_6^2  = 0 \\
\text{(E13)} = g_2 \text{(E1)} \quad A_7^2 X_2^2 - A_{8}^2 X_3^2 + A_2^2 X_4^2 - A_1^2 X_5^2 = 0 \\
\text{(E14)} = g_3 \text{(E1)} \quad A_{9}^2 X_2^2 - A_4^2 X_3^2 - A_5^2 X_4^2 + A_1^2 X_6^2 = 0 \\
\text{(E15)} = h_3 \text{(E14)} \quad A_6^2 X_2^2 - A_3^2 X_3^2 - A_5^2 X_5^2 + A_2^2 X_6^2 = 0
\end{align*}
\subsection{Twenty of eighty Rosenhain's hyperplanes}
By the action of $2$-torsion points (Remark \ref{Rem-sign}), we can get $80$ Rosenhain's hyperplanes from 
the following $20$ ones.
\begin{align*}
H1 :& \quad A_{1,3,6} X_4 + A_{2,4,9} X_5 + A_{5,7,8} X_6 =0, \quad D_{1,3,6,14} + E_{10,11,12,13}\\
H2 :& \quad A_{1,10,3} X_3 + A_{5,8,9} X_5 - A_{2,4,7} X_6 =0, \quad D_{1,3,10,13} + E_{6,11,12,14}\\
H3 :& \quad A_{2,4,10} X_3 - A_{5,6,8} X_4 - A_{1,3,7} X_6 =0, \quad D_{1,3,7,16} +E_{9,11,12,15} \\
H4 :& \quad A_{5,8,10} X_3 + A_{2,4,6} X_4 + A_{1,3,9} X_5 =0, \quad D_{1,3,9,15} +E_{7,11,12,16} \\
H5 :& \quad A_{1,6,10} X_2 + A_{4,5,7} X_5 - A_{2,8,9} X_6 =0, \quad D_{1,6,10,12} + E_{3,11,13,14} \\
H6 :& \quad A_{2,9,10} X_2 - A_{3,5,7} X_4 - A_{1,6,8} X_6 =0, \quad D_{1,6,8,16} +E_{4,11,13,15} \\
H7 :& \quad A_{5,7,10} X_2 + A_{2,3,9} X_4 + A_{1,4,6} X_5 =0, \quad D_{1,4,6,15} + E_{8,11,13,16}\\
H8 :& \quad A_{6,8,9} X_2 - A_{3,4,7} X_3 + A_{1,2,10} X_6 =0, \quad D_{1,2,10,16} + E_{5,11,14,15} \\
H9 :& \quad A_{4,6,7} X_2 - A_{3,8,9} X_3 + A_{1,5,10} X_5 =0, \quad D_{1,5,10,15} +E_{2,11,14,16} \\
H10 :& \quad A_{3,7,9} X_2 - A_{4,6,8} X_3 + A_{2,5,10} X_4 =0, \quad D_{2,5,10,14} +E_{2,5,10,14} \\
H11 :& \quad A_{3,6,10} X_1 + A_{2,7,8} X_5 - A_{4,5,9} X_6 =0, \quad D_{3,6,10,11} +E_{3,6,10,11} \\
H12 :& \quad A_{4,9,10} X_1 - A_{1,7,8} X_4 - A_{3,5,6} X_6 =0, \quad D_{3,5,6,16} +E_{4,9,10,11} \\
H13 :& \quad A_{7,8,10} X_1 + A_{1,4,9} X_4 + A_{2,3,6} X_5 =0, \quad D_{2,3,6,15} +E_{7,8,10,11} \\
H14 :& \quad A_{5,6,9} X_1 + A_{1,2,7} X_3 - A_{3,4,10} X_6 =0, \quad D_{1,2,7,13} +E_{8,12,14,15} \\
H15 :& \quad A_{2,6,7} X_1 + A_{1,5,9} X_3 - A_{3,8,10} X_5 =0, \quad D_{1,5,9,13} +E_{4,12,14,16} \\
H16 :& \quad A_{1,7,9} X_1 + A_{2,5,6} X_3 - A_{4,8,10} X_4 =0, \quad D_{2,5,6,13} +E_{4,8,10,14} \\
H17 :& \quad A_{3,4,5} X_1 + A_{1,2,8} X_2 + A_{6,9,10} X_6 =0, \quad D_{1,2,8,12} +E_{7,13,14,15} \\
H18 :& \quad A_{2,3,8} X_1 + A_{1,4,5} X_2 + A_{6,7,10} X_5 =0, \quad D_{1,4,5,12} +E_{9,13,14,16} \\
H19 :& \quad A_{1,4,8} X_1 + A_{2,3,5} X_2 + A_{7,9,10} X_4 =0, \quad D_{2,3,5,12} +E_{7,9,10,14} \\
H20 :& \quad A_{1,2,5} X_1 + A_{3,4,8} X_2 + A_{6,7,9} X_3 =0, \quad D_{1,2,5,11} +E_{10,14,15,16} \\
\end{align*}

\end{document}